\input amstex
\documentstyle {amsppt}
\UseAMSsymbols \vsize 18cm \widestnumber\key{ZZZZZ}

\catcode`\@=11
\def\displaylinesno #1{\displ@y\halign{
\hbox to\displaywidth{$\@lign\hfil\displaystyle##\hfil$}&
\llap{$##$}\crcr#1\crcr}}
\def\ldisplaylinesno #1{\displ@y\halign{
\hbox to\displaywidth{$\@lign\hfil\displaystyle##\hfil$}&
\kern-\displaywidth\rlap{$##$} \tabskip\displaywidth\crcr#1\crcr}}
\catcode`\@=12

\refstyle{A}

\font\main=cmsy10 at 10pt

\font \fin=lasy8 at 15.4 pt

\def \R{\mathop{\hbox{\main R}}\nolimits}

\topmatter
\title A note on standard modules and Vogan $L$-packets \endtitle

\rightheadtext{Standard modules and Vogan $L$-packets}
\author Volker Heiermann  \endauthor
\address Aix-Marseille Universit\'e, CNRS, Centrale Marseille, I2M, UMR 7373, 13453 Marseille, France;
\endaddress
\email volker.heiermann\@univ-amu.fr \endemail

\thanks The author has benefitted from a grant of Agence
Nationale de la Recherche with reference ANR-13-BS01-0012 FERPLAY.
\endthanks

\abstract Let $F$ be a non-Archimedean local field of characteristic $0$, let $G$ be the group of $F$-rational points of a connected reductive group defined over $F$ and  let $G'$ be the group  of $F$-rational points of its quasi-split inner form. Given standard modules $I(\tau ,\nu )$ and $I(\tau ',\nu ')$ for $G$ and $G'$ respectively with $\tau '$ a generic tempered representation, such that the Harish-Chandra $\mu $-function of a representation in the supercuspidal support of $\tau $ agrees with the one of a generic essentially square-integral representation in some Jacquet module of $\tau '$ (after a suitable identification of the underlying spaces under which $\nu =\nu '$), we show that $I(\tau ,\nu )$ is irreducible whenever $I(\tau ',\nu ')$ is.

The conditions are satisfied if the Langlands quotients $J(\tau ,\nu )$ and $J(\tau ',\nu ')$ of respectively $I(\tau ,\nu )$ and $I(\tau ',\nu ')$ lie in the same Vogan $L$-packet (whenever this Vogan $L$-packet is defined), proving that, for any Vogan $L$-packet, all the standard modules with Langlands quotient in a given Vogan $L$-packet are irreducible, if and only if this Vogan $L$-packet contains a generic representation.

This result for generic Vogan $L$-packets was proven for quasi-split orthogonal and symplectic groups by Moeglin-Waldspurger and used in their proof of the general case of the local Gan-Gross-Prasad conjectures for these groups.

\endabstract

\endtopmatter
\document
Let $F$ be a non-Archimedean local field of characteristic $0$. Denote by $G$ the group of $F$-rational points of a connected reductive group defined over $F$.

Fix a minimal $F$-parabolic subgroup $P_0=M_0U_0$ of $G$. A $F$-parabolic subgroup $P=MU$ will be called standard, if it contains $P_0$ and $M$ contains $M_0$ and semi-standard, if only $M$ contains $M_0$. One says then respectively that $M$ is a standard or semi-standard $F$-Levi subgroup. If $M$ is any semi-standard $F$-Levi subgroup of $G$, a standard $F$-parabolic subgroup of $M$ will be a $F$-parabolic subgroup of $M$ which contains $P_0\cap M$ and analog for the $F$-Levi subgroups of $M_0$.

Let $P=MU$ be a standard $F$-parabolic subgroup of $G$ and $T_M$ the maximal split torus in the center of $M$. We will write $a_M^*$ for the dual of the real Lie-algebra $a_M$ of $T_M$, $a_{M,\Bbb C}^*$ for its complexification and $a_M^{*+}$ for the positive Weyl chamber in $a_M^*$ defined with respect to $P$. The set of reduced roots for the action of $T_M$ in the Lie algebra of $G$ will be denoted $\Sigma (T_M)$, the subset of positive roots $\Sigma (P)$, and the corresponding sets of coroots respectively $\Sigma ^{\vee }(T_M)$ and $\Sigma^{\vee }(P)$. They will be identified with subsets of $a_M^*$ and of $a_M$ respectively.

Following \cite{W}, one defines a map $H_M:M\rightarrow a_M$, such that, for every $F$-rational character $\chi\in a_M^*$ of $M$, $\vert\chi (m)\vert _F=q^{-\langle\chi ,H_M(m)\rangle}$. If $\pi $ is a smooth representation of $M$ and $\nu \in a_{M,\Bbb C}^*$, we denote by $\pi _{\nu }$ the smooth representation of $M$ defined by $\pi _{\nu }(m)=q^{-\langle\nu ,H_M(m)\rangle}\pi (m)$. (Remark that, although the sign in the definition of $H_M$ has been changed compared to the one due to Harish-Chandra, the meaning of $\pi _{\nu }$ is unchanged.) The symbol $i_P^G$ will denote the functor of parabolic induction normalized such that it sends unitary representations to unitary representations, $G$ acting on its space by right translations.

Let $\tau $ be an irreducible tempered representation of $M$ and $\nu\in a_M^{*+}$. Then the induced representation $i_P^G\tau _{\nu }$ is called a standard module and will be denoted $I(P,\tau ,\nu )$ in the sequel. It has a unique irreducible quotient $J(P,\tau ,\nu )$, the so-called Langlands quotient. By the Langlands classification theorem, every irreducible smooth representation $\pi $ of $G$ appears as a Langlands quotient of a unique standard module $I(P,\tau ,\nu )$. One says then that $(P,\tau ,\nu )$ are Langlands data for $\pi $ and writes $(P_{\pi },\tau _{\pi },\nu _{\pi })$ if one wants to underline the link to $\pi $.

Let $G'$ be the group of $F$-rational points of a quasi-split connected reductive group defined over $F$. Fix a $F$-Borel subgroup $B'=T'V'$, which defines the notions of standard $F$-parabolic, standard $F$-Levi, and semi-standard $F$-Levi subgroup of $G'$ and of its semi-standard $F$-Levi subgroups. Suppose that there is a standard $F$-parabolic subgroup $P_0'=M_0'U_0'$ of $G'$ so that $(a_{M'_0}^*,\Sigma (T_{M'_0}), a_{M_0'}, \Sigma^{\vee }(T_{M'_0}))$ is a root datum and that the based root data $(a_{M'_0}^*,\Sigma (T_{M'_0}), a_{M_0'}, \Sigma^{\vee }(T_{M'_0}),\Delta (P_0'))$ and $(a_{M_0}^*,\Sigma (T_{M_0}), a_{M_0}, \Sigma^{\vee }(T_{M_0}),\Delta (P_0))$ are isomorphic. The isomorphism is a pair $(\varphi ,\varphi ^{\vee })$, where $\varphi $ is an isomorphism between $a_{M_0}^*$ and $a_{M'_0}^*$ sending $\Sigma (T_{M_0})$ onto $\Sigma (T_{M'_0})$ and $\Delta (P_0)$ onto $\Delta (P_0')$, and $\varphi ^{\vee }$ is the inverse dual map that sends $\Sigma ^{\vee }(T_{M_0})$ onto $\Sigma ^{\vee }(T_{M'_0})$. (One may think of $G'$ as an inner form of $G$ (and it is then unique), but this is not required below.)

In particular, the isomorphism $\varphi $ of based root data defines a bijection between standard $F$-parabolic subgroups $P'=M'U'$ of $G'$, $P'\supseteq P_0'$, and standard $F$-parabolic subgroups $P=MU$ of $G$ with an isomorphism $a_M^*\rightarrow a_{M'}^*$, also denoted $\varphi $ in the sequel, that sends the set of positive roots $\Sigma (P)$ onto $\Sigma (P')$. We will say that the standard parabolic subgroups $P=MU$ and $P'=M'U'$ \it correspond to each other by $\varphi $. \rm Remark that the inverse dual map $\varphi ^\vee:a_{M}\rightarrow a_{M'}$ sends then $\Sigma^{\vee }(P)$ onto $\Sigma^{\vee }(P')$.

Denote by $W'$ the Weyl group of $G'$ defined with respect to the maximal $F$-split torus $T'$ in $G'$ and by ${w'}_0^G$ the longest element in $W'$. After changing the splitting in $V'$, for any generic representation $\pi '$  of $G'$, one can always find a non degenerate character $\psi '$ of $V'$, which is compatible with ${w'}_0^G$, such that $\pi '$ is $\psi '$-generic \cite{Sh, section 3}. For any semi-standard Levi-subgroup $M'$ of $G'$, we will still denote by $\psi '$ the restriction of $\psi '$ to $M'\cap V'$. It is compatible with ${w'}_0^M$. If we write in the sequel that a representation of a $F$-semi-standard Levi subgroup of $G'$ is $\psi '$-generic, then we always mean that $\psi '$ is a non degenerate character of $V'$ with the above properties.

Fix for the rest of this introduction a standard $F$-parabolic subgroup $P=MU$ of $G$ and a standard $F$-parabolic subgroup $P'=M'U'$ of $G'$ that correspond to each other by the above isomorphism of based root data $\varphi $.

Let $\tau $ be an irreducible tempered representation of $M$, $\sigma $ a unitary irreducible supercuspidal representation of a standard $F$-Levi subgroup $M_1$ of $M$ and $\lambda _{\tau }\in a_{M_1}^{M*}$, so that $\sigma _{\lambda _{\tau }}$ lies in the supercuspidal support of $\tau $. Denote by $P_1'=M_1'U_1'$ the standard $F$-parabolic subgroup of $G'$ that corresponds to $P_1$ by the isomorphism of root data $\varphi $. Assume that there is a $\psi '$-generic discrete series representation $\sigma '$ of $M'_1$, such that Harish-Chandra's $\mu $-functions $\mu ^M$ and $\mu ^{M'}$ become equal under $\varphi $ on the inertial orbits of $\sigma $ and $\sigma '$ , i.e. $\mu ^M(\sigma _{\lambda })=\mu ^{M'}(\sigma '_{\varphi(\lambda )})$ for all $\lambda\in a_{M_1}^*$.

Our first result is then (cf. lemma {\bf 2.5}) that $i_{P_1'\cap M'}^{M'}\sigma '_{\varphi(\lambda _{\tau })}$ has a $\psi '$-generic tempered sub-quotient $\tau '$.

The aim of this note is further to show that under these conditions the standard module $I(P,\tau ,\nu )$ is irreducible, whenever the standard module $I(P',\tau' ,\varphi(\nu ))$ is irreducible (cf. theorem {\bf 2.7}). We show in the annex {\bf A.} that the above conditions are satisfied if $G'$ is a quasi-split classical group and $\pi $ and $\pi '$ are two members of a generic Vogan $L$-packet, $\pi '$ being the generic member, and $(P,\tau ,\nu )$ and $(P',\tau ',\nu ')$ are Langlands data for them. (For a general quasi-split group $G'$, one may add this property to the conjectural properties of Vogan $L$-packets.) It follows, using the standard modules conjecture for $G'$ proven in its whole generality in \cite{HM} (subsequel to \cite{HO}) that $I(P,\tau ,\nu )$ is always irreducible when the Langlands quotient $J(P,\tau ,\nu )$ lies in a $\psi '$-generic Vogan $L$-packet (cf. theorem {\bf 3.2}).

The above result for generic Vogan $L$-packets was proven for quasi-split orthogonal and symplectic groups by Moeglin-Waldspurger in \cite{MW, 2.12} and used in their proof of the general case of the local Gan-Gross-Prasad conjectures for these groups. (The restrictions that we make in {\bf 3.1} on the literature hold also for their work.)

As the Vogan $L$-packets for unitary groups (with their conjectured properties) are known following the work of Mok \cite{Mk} and Kaletha-Minguez-Shin-White \cite{KMSW}, our results apply in particular to standard modules corresponding to generic Vogan $L$-packets of unitary groups. This should be useful if one wants to prove the general case of the local Gan-Gross-Prasad conjectures for unitary groups, the conjectures for tempered representations having been established by Beuzart-Plessis \cite{BP}.

The method of proof of the present note is a generalization of the proof of the standard modules conjecture in \cite{HM}, based on the basic principle in \cite{Mu}, following the observation that everything is determined by the geometric combinatorial configuration of the parameters with respect to Harish-Chandra's $\mu $-function. Remark that the result in \cite{HM} depended on the tempered $L$-function conjecture which has been proven in its whole generality in \cite{HO}. (The fact that the theory of $L$-functions is only established in characteristic $0$ is the reason why we had to make this assumption on the characteristic of $F$. Thanks to the work of L. Lomel\'\i\  and others, this obstacle may be lifted soon.)

\null{\bf 1.} Recall that $G'$ denotes the set of $F$-rational points of a quasi-split connected reductive group defined over $F$ and that one has fixed a Borel subgroup $B'=T'V'$ of $G'$ and a non-degenerate character $\psi '$ of $V'$.

\null{\bf 1.1 Theorem:} (Rodier) \it Let $P'=M'U'$ be a standard $F$-parabolic subgroup of $G'$ and $\tau '$ an irreducible smooth representation of $M'$.

The representation $i_{P'}^{G'}\tau '$ has an irreducible $\psi '$-generic subquotient, if and only if $\tau '$ is $\psi '$-generic. This irreducible $\psi '$-generic subquotient is then unique. \rm

\null {\bf 1.2} Let $M'$ be a semi-standard Levi subgroup of $G'$ and $\tau '$ an irreducible generic tempered representation of $M'$. In \cite{Sh3}, F. Shahidi attached to each irreducible component $r_i'$ of the adjoint action of the $L$-group $^LM'$ of $M'$ on $Lie(\ ^LU')$ an $L$-function $L(s,\tau ',r_i')$, which is a meromorphic function in the complex variable $s$.

\null {\bf Theorem:} \cite{HO, 5.1} \it The meromorphic function $L(s,\tau ',r_i')$ is regular for $\R (s)>0$. \rm

\null {\bf 1.3} The following statement was proven in \cite{HM, 2.2} modulo theorem {\bf 1.2} (which was proven in its whole generality later).

\null {\bf Proposition:} \it Let $\sigma '$ be a $\psi '$-generic supercuspidal representation of the standard Levi factor of a standard parabolic subgroup $P'$ of $G'$.

Suppose that $i_{P'}^{G'}\sigma '$ has an irreducible subquotient that is a tempered representation. Then, the $\psi '$-generic sub-quotient of $i_{P'}^{G'}\sigma '$ is a tempered representation.

\null{\bf 1.4 Corollary:} \it Let $\sigma '$ be a $\psi '$-generic representation of the standard Levi factor of a standard parabolic subgroup $P'$ of $G'$.

Suppose that $i_{P'}^{G'}\sigma '$ has an irreducible sub-quotient that is a tempered representation. Then, the generic subquotient of $i_{P'}^{G'}\sigma '$ is a tempered representation.

\null Proof: \rm Fix a standard parabolic subgroup $P_1'=M_1'U_1'$ of $G'$ with $M_1'\subseteq M'$ and an irreducible supercuspidal representation $\sigma _1'$ of $M_1'$, such that $\sigma '$ is a sub-representation of $i_{P_1'\cap M'}^{M'}\sigma _1'$. The representation $\sigma _1'$ is $\psi '$-generic by {\bf 1.1}. An irreducible sub-quotient of $i_{P'}^{G'}\sigma '$ is necessarily a sub-quotient of $i_{P_1'}^{G'}\sigma _1'$. In particular, $i_{P_1'}^{G'}\sigma _1'$ has a tempered sub-quotient. By {\bf 1.3}, the $\psi '$-generic sub-quotient $\pi $ is tempered, too. By uniqueness of the $\psi '$-generic sub-quotient, $\pi $ is also the $\psi '$-generic subquotient of $i_{P'}^{G'}\sigma '$. This proves the corollary. \hfill{\fin 2}

\null{\bf 2.} Recall that $G$ denotes the $F$-rational points of a connected reductive group defined over $F$ and that one has fixed a minimal $F$-parabolic subgroup $P_0=M_0U_0$ of $G$ and a $F$-parabolic subgroup $P'_0=M'_0U'_0$ of the quasi-split group $G'$ so that the corresponding based root spaces are isomorphic by an isomorphism $\varphi $ (see introduction).

\null
{\bf 2.1} The Harish-Chandra $\mu $-function is the main ingredient of the Plancherel density for a $p$-adic reductive group $H$ \cite{W}. It assigns to every square-integrable representation of a Levi subgroup a complex number and can be analytically extended to a meromorphic function on the space of essentially square-integrable representations of Levi subgroups. What is important here, is the following: if $Q=NV$ is a parabolic subgroup of a connected reductive group $H$ over $F$ and $\sigma $ an irreducible unitary supercuspidal representation of $N$, then the Harish-Chandra's $\mu $-function $\mu ^H$ corresponding to $H$ defines a meromorphic function $a_{N,\Bbb C}^*\rightarrow\Bbb C$, $\lambda\mapsto\mu^H(\sigma _{\lambda})$, on the $\Bbb C$-vector space $a_{N,\Bbb C}^*$ which can be written in the form $$\leqno{(2.1)}\qquad\displaystyle{\mu^H(\sigma _\lambda )=f(\lambda )\prod _{\alpha\in\Sigma(Q)}{(1-q^{\langle\alpha^{\vee },\lambda \rangle})(1-q^{-\langle\alpha^{\vee },\lambda \rangle})\over (1-q^{\epsilon_{\alpha ^{\vee }}+\langle\alpha^{\vee },\lambda \rangle})(1-q^{\epsilon_{\alpha ^{\vee }}-\langle\alpha^{\vee },\lambda \rangle})}}\ \ \hbox{\rm (cf. \cite{H1, Si})},$$ where $f$ is a meromorphic function without poles and zeroes on $a_N^*$ and the $\epsilon _{\alpha ^{\vee }}$ are non negative rational numbers with $\epsilon _{\alpha ^{\vee }}=\epsilon _{{\alpha '}^{\vee }}$ if $\alpha $ and $\alpha '$ are conjugate.

\null{\bf 2.2} \it Definition: \rm A point $\sigma _{\lambda _0}$, $\lambda _0\in a_N^*$, is called a residue point for $\mu ^H$, if $$\vert\{\alpha\in\Sigma (Q)\vert \langle\alpha^{\vee },\lambda _0\rangle =\pm\epsilon _{\alpha ^{\vee }}\}\vert - 2\vert\{\alpha\in\Sigma (Q)\vert \langle\alpha^{\vee },\lambda_0\rangle =0\}\vert =dim(a_N^*/a_G^*).$$

\null{\bf 2.3 Theorem:} \cite{H1} \it The induced representation $i_Q^H\sigma_{\lambda _0}$ has a sub-quotient which is square integrable, if and only if $\sigma _{\lambda _0}$ is a residue point of $\mu ^H$ and the projection of $\lambda _0$ on $a_H^*$ is $0$. \rm

\null{\bf 2.4 Lemma:} \it Let $\sigma $ be a unitary irreducible supercuspidal representation of the standard $F$-Levi subgroup $M$ of a standard parabolic subgroup $P$ of $G$. Denote by $P'=M'U'$ the $F$-standard parabolic subgroup of $G'$ that corresponds to $P$ by the isomorphism of basic root data $\varphi $.

Assume that there is a $\psi '$-generic discrete series representation $\sigma '$ of $M'$ such  that the Harish-Chandra $\mu $-functions of $\sigma $ and $\sigma '$ are directly proportional under $\varphi $,  i.e. that there is some complex number $c\in\Bbb C$ such that $\mu ^G(\sigma _{\lambda })=c\mu ^{G'}(\sigma '_{\varphi (\lambda )})$ for all $\lambda\in a_{M,\Bbb C}^*$.

Let $\lambda _0\in a_M^*$. If $\sigma _{\lambda _0}$ is a residue point for $\mu ^G$, then $\sigma '_{\varphi (\lambda _0)}$ is a residue point for $\mu ^{G'}$ and vice versa. In addition, the projection of $\lambda _0$ onto $a_G^*$ is zero if only if the projection of $\varphi (\lambda _0)$ is zero on $a_{G'}^*$.

\null Remark: \rm Although the notion of a residue point is only defined for $\mu $-functions with respect to supercuspidal representations, the statement makes sense because the $\mu $-function with respect to $\sigma '_{\lambda '}$, $\lambda '\in a_{M'}^*$, behaves by our assumptions like the one of a supercuspidal representation.

\null\it Proof: \rm By the formula (2.1),  there are a meromorphic function $f$ without poles and zeroes on $a_N^*$ and non negative rational numbers $\epsilon _{\alpha ^{\vee }}$, $\epsilon _{\alpha ^{\vee }}=\epsilon _{{\alpha '}^{\vee }}$ if $\alpha $ and $\alpha '$ are conjugate, such that, for $\lambda'\in a_{M'}^*$, $$\mu^G(\sigma _{\varphi^{-1}(\lambda ')})=f(\varphi^{-1}(\lambda '))\prod _{\alpha\in\Sigma(P)}{(1-q^{\langle\alpha^{\vee },\varphi^{-1}(\lambda ')\rangle})(1-q^{-\langle\alpha^{\vee },\varphi^{-1}(\lambda ')\rangle})\over (1-q^{\epsilon_{\alpha ^{\vee }}+\langle\alpha^{\vee },\varphi^{-1}(\lambda ')\rangle})(1-q^{\epsilon_{\alpha ^{\vee }}-\langle\alpha^{\vee },\varphi^{-1}(\lambda ')\rangle})}.$$ Applying the dual of $\varphi ^{-1}$, this becomes $$f(\varphi^{-1}(\lambda '))\prod _{\alpha\in\Sigma(P)}{(1-q^{\langle\varphi ^{\vee }(\alpha^{\vee }),\lambda '\rangle})(1-q^{-\langle\varphi ^{\vee }(\alpha^{\vee }),\lambda '\rangle})\over (1-q^{\epsilon_{\alpha ^{\vee }}+\langle\varphi ^{\vee }(\alpha^{\vee }),\lambda '\rangle})(1-q^{\epsilon_{\alpha ^{\vee }}-\langle\varphi ^{\vee }(\alpha^{\vee }),\lambda' \rangle})}.$$  By our assumptions, this equals $\mu^{G'}(\sigma '_{\lambda '})$ up to a constant factor. Comparing this expression with the expression for $\mu^{G'}(\sigma '_{\lambda '})$ given by (2.1), one deduces that $\epsilon_{\alpha ^{\vee }}=\epsilon_{\varphi^{\vee }(\alpha ^{\vee })}$. Looking at the zero and polar hyperplanes in both cases, it becomes clear that $\sigma'_{\lambda '}$ is a residue point of $\mu ^{G'}$ if and only if $\sigma _{\varphi^{-1}(\lambda ')}$ is one for $\mu ^G$.

The projection of an element $\lambda \in a_M^*$ onto $a_G^*$ is $0$, if and only if $\lambda $ is a linear combination of elements in $\Sigma (T_M)$. As $\varphi $ induces a bijection between $\Sigma (T_M)$ and $\Sigma (T_{M'})$, this is equivalent to say that $\varphi (\lambda )$ is a linear combination of elements in $\Sigma (T_{M'})$. Consequently, the projection of $\lambda $ onto $a_G^*$ is $0$, if and only if the projection of $\varphi (\lambda )$ onto $a_{G'}^*$ is $0$. \hfill{\fin 2}

\null{\bf 2.5 Lemma:} (with the notations and assumptions of {\bf 2.4}) \it Let $\lambda\in a_M^*$ and assume that $i_P^G\sigma _{\lambda }$ has a tempered sub-quotient. Then, $i_{P'}^{G'}\sigma '_{\varphi(\lambda )}$ has a $\psi '$-generic tempered sub-quotient.

\null Proof: \rm Let $\tau $ be a tempered sub-quotient of $i_P^G\sigma _{\lambda }$. Then, it is a sub-representation of a representation induced by a discrete series representation $\tau _2$ of some standard $F$-Levi subgroup $M_2$ of $G$ and $\sigma _{\lambda }$ is conjugated to an element in the supercuspidal support of $\tau _2$. After conjugating $\sigma _{\lambda }$ and $M$ by a Weyl group element (which does not effect the Jordan-H\"older sequence for $i_P^G\sigma _{\lambda }$ neither for $i_{P'}^{G'}\sigma '_{\varphi (\lambda )}$ - if $\sigma _{\lambda }$ is conjugated by $w$, $\sigma '_{\varphi (\lambda )}$ has to be conjugated by $w':=\varphi\circ w \circ \varphi ^{-1}$), we can assume that $M\subseteq M_2$ and consequently $\sigma _{\lambda }$ lies in the supercuspidal support of $\tau _2$. By {\bf 2.3}, this implies that $\sigma _{\lambda }$ is a residue point of the Harish-Chandra $\mu $-function $\mu ^{M_2}$ and that the projection of $\Re(\lambda )$ onto $a_{M_2}^*$ is trivial.

Denote by $M_2'$ the standard $F$-Levi subgroup of $G'$ that corresponds to $M_2$ by the isomorphism of root data $\varphi $. By {\bf 2.4}, $\sigma '_{\varphi(\lambda )}$ is a residue point of $\mu ^{M_2'}$ and the projection of $\Re(\lambda )$ on $a_{M_2'}^*$ is trivial.

It remains to show that $i_{P'\cap M_2'}^{M_2'}\sigma '_{\varphi(\lambda )}$ has a discrete series sub-quotient. If $\sigma '$ is not a supercuspidal representation, one cannot immediately apply theorem {\bf 2.3}. Let $P_1'=M_1'U_1'$ be a standard parabolic subgroup of $G'$ contained in $P'$, $\sigma _1'$ a $\psi '$-generic unitary supercuspidal representation of $M_1'$ and $\lambda _2'\in a_{M_1'}^{M'*}$ such that $\sigma _{1,\lambda _2'}'$ lies in the supercuspidal support of $\sigma '$. As $\sigma '$ is a square-integrable representation, $\sigma _{1,\lambda _2'}$ is a residue point of $\mu ^{M'}$. One remarks that there is a complex number $c$, such that $c\mu ^{M_2'}(\sigma '_{\lambda '} )=(\mu^{M_2'}/\mu ^{M'})(\sigma '_{1,\lambda _2'+\lambda '})$  as functions in $\lambda '\in a_{M'}^*$ (cf. lemma {\bf A.3} in the attachment). As $\sigma'_{\varphi(\lambda )}$ is a residue point for $\mu ^{M_2'}$, one sees that $\sigma '_{1,\lambda _2'+\varphi(\lambda )}$ is also a residue point for $\mu ^{M_2'}$ and the projection on $a_{M_2'}^*$ is zero. From {\bf 2.3}, it follows that $i_{P_1'\cap M_2'}^{M_2'}\sigma '_{1,\lambda _2'+\varphi(\lambda )}$ has a square-integrable sub-quotient and consequently $i_{P_1'}^{G'}\sigma '_{1,\lambda _2'+\varphi(\lambda )}$ has a tempered sub-quotient. As $\sigma '_{1,\lambda _2'+\varphi(\lambda )}$ is $\psi '$-generic, $i_{P_1'}^{G'}\sigma '_{1,\lambda _2'+\varphi(\lambda )}$ has by {\bf 1.4} a $\psi '$-generic tempered sub-quotient, which is by its uniqueness the $\psi '$-generic sub-quotient of $i_{P'}^{G'}\sigma '_{\varphi(\lambda )}$. This proves the lemma.
\hfill{\fin 2}

\null{\bf 2.6}\ \it D\'efinition: \rm Let $I(P,\tau ,\nu )$ be a standard module for $G$ and $I(P',\tau ',\nu ')$ a $\psi '$-generic standard module for $G'$. We will say that \it $I(P,\tau ,\nu )$ is related to $I(P',\tau ',\nu ')$ (by $\varphi $), \rm if $\nu '=\varphi(\nu )$, the two $F$-standard parabolic subgroups $P=MU$ and $P'=M'U'$ correspond to each other by $\varphi $, and there are in addition standard parabolic subgroups $P_1=M_1U_1$ and $P_1'=M_1'U_1'$ of $G$ and $G'$ respectively, a unitary supercuspidal representation $\sigma $ of $M_1$, a $\psi '$-generic discrete series representation $\sigma '$ of $M_1'$ such that $P_1$ and $P_1'$ correspond to each other by $\varphi $, $\mu ^M(\sigma _{\lambda })$ and $\mu ^{M'}(\sigma '_{\varphi(\lambda )})$ are directly proportional as functions of $\lambda\in a_{M_1}^*$ (i.e. there is some complex number $c$ such that $\mu ^M(\sigma _{\cdot })=c\mu ^{M'}(\sigma '_{\varphi(\cdot )})$), and there is $\lambda _1\in a_{M_1}^*$ such that $\tau $ is a subquotient of $i_{P_1\cap M}^{M}\sigma _{\lambda _1}$ and $\tau '$ a subquotient of $i_{P_1'\cap M'}^{M'}\sigma '_{\varphi (\lambda _1)}$.

\null\it Remark: \rm 1) In other words, the supercuspidal representation $\sigma _{\lambda _1}$ of $M_1$ is an element of the supercuspidal support of $\tau $ and $\sigma '$ is a $\psi '$-generic discrete series representation of the standard Levi subgroup $M_1'$ of $G'$ that corresponds to $M_1'$ and with $\mu $-function directly proportional to the one of $\sigma $ via the identification $\varphi $.

2) There is a distinguished class of inner forms of $G'$, called pure inner forms, defined by Vogan \cite{V}. It is expected that $I(P,\tau ,\nu )$ is related to $I(P',\tau ',\nu ')$, if $G'$ is a pure inner form of $G$ and $J(P,\tau ,\nu )$ and $J(P',\tau ',\nu ')$ lie in the same Vogan $L$-packet (defined in \cite{V}). If $G=G'$, this means that these representations lie in the same (Langlands) $L$-packet. In addition, if $G'$ is a pure inner form of $G$, then every standard module $I(P,\tau ,\nu )$ for $G$ should be related to a $\psi '$-generic standard module for $G'$. All this is known now for quasi-split classical groups (orthogonal, symplectic or unitary) by the stabilization of the Arthur-Selberg trace formula and work of different authors on it, see {\bf 3.1} for more details.

By a generalization of the Jacquet-Langlands correspondence, similar properties should hold for a general inner form $G$ of $G'$.

\null{\bf 2.7 Theorem:} \it Let $I(P',\tau ',\nu ')$ be a $\psi '$-generic standard module for $G'$. If $I(P',\tau ',\nu ')$ is irreducible, then every standard module $I(P,\tau ,\nu )$ of $G$ that is related to $I(P',\tau ',\nu ')$ is also irreducible. In particular, if $J(P',\tau ',\nu ')$ is $\psi'$-generic, each standard module $I(P,\tau ,\nu )$ of $G$ that is related to $I(P',\tau ',\nu ')$ is irreducible.

\null \it Proof: \rm  Fix data $P_1=M_1U_1,\sigma _1, \lambda _1$ and $P_1'=M_1'U_1',\sigma _1'$ associated to  $I(P,\tau, \nu)$ and $I(P',\tau ',\nu ')$ as in definition {\bf 2.6}. In particular, $\sigma_{\lambda _1}$ lies in the supercuspidal support of $\tau $ and $\sigma'_{\varphi(\lambda _1')}$ in the supercuspidal support of $\tau '$.

If $I(P,\tau ,\nu )$ is reducible, then there is an irreducible subquotient $\pi _2$ of $I(P,\tau ,\nu )$ that is not isomorphic to $J(P,\tau ,\nu )$. By \cite{BW, chapter XI, 2.13}, there is a Langlands data $(P_2,\tau _2,\nu _2)$ with $\nu _2<\nu $ in $a_T^*$ (the precise definition of this order does not matter here, only that it is stable by the isomorphism of based root data $\varphi $), so that $\pi _2=J(P_2,\tau _2,\nu _2)$. Remark that $\sigma _{\lambda_1+\nu }$ lies in the supercuspidal support of $\pi $, $\pi _2$ and $\tau _{\nu }$, and that it is at least conjugate to an element in the supercuspidal support of $\tau _{2,\nu _2}$. This element in the supercuspidal support of  $\tau _{2,\nu _2}$ can be written in the form $(w\sigma )_{w\lambda _1+w\nu }$ where $w$ is some Weyl group element which acts by conjugation so that $wM_1:=wM_1w^{-1}$ is a standard $F$-Levi subgroup. Put $w'=\varphi\circ w\circ\varphi ^{-1}$. Then, the standard $F$-Levi subgroup of $G'$ that corresponds to $wM_1$ by the isomorphism of root data $G$ is $w'M_1'$.

As $\sigma $ and $\sigma '$ have directly proportional $\mu $-functions by $\varphi $, the representation $i_{M'\cap P_1'}^{M'}\sigma '_{\varphi(\lambda _1)}$ has by lemma {\bf 2.5} a $\psi '$-generic tempered sub-quotient $\tau '$. As the standard $F$-Levi subgroups $wM_1$ and $w'M_1'$ correspond to each other by the isomorphism of root  data $\varphi $ and Plancherel measures are invariant by conjugation, $w\sigma $ and $w'\sigma '$ have also directly proportional $\mu $-functions by $\varphi $. Denote by $P_2'=M_2'U_2'$ the standard $F$-parabolic subgroup of $G'$ that corresponds to $P_2=M_2U_2$ and by $P_1'$ the standard $F$-parabolic subgroup of $G'$ with Levi subgroup $M_1'$. It follows from lemma {\bf 2.5} that the induced representation $i_{wP_1'\cap M_2'}^{M_2'}\sigma '_{\varphi (w\lambda _1+w\nu -\nu _2)}$ has a $\psi '$-generic tempered sub-quotient $\tau _2'$. Consequently, the $\psi '$-generic sub-quotient $\pi ^{gen}$ of $I(P',\tau ',\varphi(\nu ))$ is also a $\psi '$-generic sub-quotient of $I(P_2',\tau _2',\varphi(\nu _2))$. As $\nu _2<\nu $ implies $\varphi (\nu _2)<\varphi (\nu )$, this proves that $\nu _{\pi ^{gen}}<\varphi(\nu )$, showing that $\pi ^{gen}$ cannot be the Langlands quotient of $I(P',\tau ',\varphi(\nu ))$. This implies that $I(P',\tau ',\varphi(\nu ))$ is reducible.

The last sentence follows from the standard modules conjecture for $G'$. \hfill{\fin 2}

\null\it Remark: \rm The opposite direction of the statement in the theorem is wrong. There are many examples of standard modules coming from the same $L$-packet, some of them being irreducible and others being reducible (see for example  \cite{T, 5.1 (vi)}).

\null {\bf 3.} We assume now that $G'$ is the set of $F$-rational points of the quasi-split inner form of $G$ and that $G$ is a pure inner form of $G'$ in the sense of Vogan \cite{V}. There exists then a natural isomorphism of based root data as mentioned in the introduction that we will still denote $\varphi $.

\null {\bf 3.1} The notion of Vogan $L$-packets is defined in \cite{V}. Remark first that $G'$ and $G$ have as inner forms of each other the same $L$-group, denoted here $\hbox {\main G}$. Two representations $\pi $ and $\pi '$ of $G$ and $G'$ respectively lie in the same Vogan $L$-packet if they correspond to the same, in general hypothetical, Langlands parameter $\eta:W_F\times SL_2(\Bbb C)\rightarrow\hbox{\main G}$. (For quasi-split classical groups, Vogan $L$-packets are known thanks to the stabilization of the trace formula in \cite{A, Mk, KMSW} and its application \cite{M1, M2}, see also \cite{H4} for some summary, only the analogue of \cite{A} for pure inner forms of orthogonal groups and of \cite{M2} for pure inner forms of quasi-split unitary groups remains to be published.)

By the expected properties of the local Langlands correspondence (see for example \cite{M1, M2} for the known case of quasi-split classical groups and \cite{H2, H3} for the general philosophy), this implies in particular that, if $(P,\tau ,\nu )$ and $(P',\tau ',\nu ')$ are Langlands data for $\pi $ and $\pi '$, then the standard parabolic subgroups $P$ and $P'$ correspond to each other, the tempered representations $\tau $ and $\tau '$ lie in the same Vogan $L$-packet, $\tau '$ will be $\psi '$-generic if $\pi '$ is $\psi '$-generic, and equality $\nu'=\varphi (\nu )$ holds.

If $G'$ is a quasi-split classical group, the following result is proven in the annex (cf. theorem {\bf A.1}):

\null{\bf Theorem:} \it Suppose that $G'$ is a quasi-split classical group. If $\sigma $ and $\sigma '$ are discrete series representations of standard Levi subgroups $M$ and $M'$ of $G$ and $G'$ respectively that correspond to each other, then their Harish-Chandra $\mu $-functions are directly proportional, i.e. there exists a complex number $c$ such that $$\mu (\sigma _{\lambda })=c\mu (\sigma '_{\varphi(\lambda )})$$ as functions in $\lambda\in a_{M,\Bbb C}^*$. \rm

\null For a general quasi-split group $G'$, one may put this theorem as a conjecture along the existence of Vogan $L$-packets. Remark that, for Langlands $L$-packets, this conjecture appeared in \cite{Sh, 9.3}. In general, some evidence can be seen from the formal degree conjecture for Vogan $L$-packets by Hiraga-Ichino-Ikeda \cite{HII}. More philosophically, one can say that $\mu $-functions are essentially defined by certain $L$-functions attached to them and which should not differ in an $L$-packet (see \cite{GT} for the case of $GSp_4$, where a more precise result is given).

\null{\bf 3.2 Theorem:} \it Assume either that $G'$ is a quasi-split classical group or that the conclusions of theorem {\bf 3.1} and the expected properties of Vogan $L$-packets hold for $G'$.

Let $\Pi $ be the Vogan $L$-packet of an irreducible smooth representation of $G'$. Fix  an $F$-parabolic subgroup $P'=M'U'$ of $G'$, an irreducible $\psi '$-generic tempered representation $\tau '$ of $M'$ and an element $\nu '$ in $a_{M'}^{*+}$ so that the Langlands quotient $J(P',\tau ',\nu ')$ lies in $\Pi $.

The following are equivalent:

(i) the standard module $I(P',\tau ',\nu ')$ is irreducible;

(ii) for every $\pi $ in $\Pi $ and Langlands data $(P,\tau ,\nu )$ for $\pi $, the standard module $I(P,\tau ,\nu )$ is irreducible.

(iii) the Langlands quotient $J(P',\tau ',\nu ')$ is $\psi '$-generic.

(iv) the Vogan $L$-packet $\Pi $ contains a generic representation.

\null Proof: \rm The existence of the $\psi '$-generic tempered representation $\tau '$ of a standard $F$-Levi subgroup $M'$ comes from the conjecture that tempered $L$-packets are generic (one can then always find a $\psi '$-generic representation \cite{Sh, 3.}). This conjecture is known for quasi-split groups thanks to the above mentioned works on the trace formula. For general $G'$, they take part of our assumptions on the Vogan $L$-packets.

As the conclusions of theorem {\bf 3.1} hold by assumption, it follows from the discussion in {\bf 3.1} and the properties of the local Langlands correspondence (known for quasi-split groups) that the standard module $I(P,\tau ,\nu )$ is related to $I(P',\tau ',\nu ')$ in the sense of definition {\bf 2.6} so that one can conclude from theorem {\bf 2.7}. \hfill{\fin 2}

\null{\bf 3.3} One can reformulate things in the following way:

\null{\bf  Theorem:} (with the assumptions of {\bf 3.2}) \it Denote by $\hbox{\main T}$ the Vogan $L$-packet of a $\psi '$-generic tempered representation $\tau '$ of the standard $F$-Levi subgroup of a standard $F$-parabolic subgroup $P'$ of $G'$. The following are equivalent:

(i) the standard module $I(P',\tau ',\nu ')$ is irreducible;

(ii) for every inner form $G$ of $G'$ and Langlands data $(P,\tau ,\nu )$ with $\tau $ in $\hbox{\main T}$ and $\nu $ corresponding to $\nu '$ by an appropriate isomorphism of based root data for $G$ and $G'$, the standard module $I(P,\tau ,\nu )$ is irreducible.

(iii) the Langlands quotient $J(P',\tau ',\nu ')$ is $\psi '$-generic. \rm

\null
\it Remark: \rm Recall that there exist many examples of standard modules which are irreducible, while the corresponding $\psi '$-generic standard module is reducible (cf. for example \cite{T, 5.1 (vi)}).

\null{\bf A. Annex}

\null The aim of this annex is to give a proof of theorem {\bf 3.1}:

\null{\bf A.1 Theorem:} \it Suppose that $G'$ is a quasi-split classical group. If $\tau ' $ is an irreducible discrete series representation of a Levi subgroup of $M'$ and $\tau $ is an element of the Vogan $L$-packet of $\tau '$, then the $\mu $-functions associated to $\tau $ and $\tau '$ become directly proportional  after identification of the real Lie algebras of the split centers of the corresponding Levi subgroups. \rm

\null The proof of this theorem will occupy the rest of this section. It will be done in several steps.

\null {\bf A.2} That the $\mu $-function is constant on $L$-packets of discrete series representations of a standard Levi subgroup of a quasi-split classical group is in \cite{Co, 7.8}. (The proof holds for every group where $L$-packets of discrete series representations are known and a linear combination of their characters gives a stable distribution - this is the case for quasi-split orthogonal and symplectic groups thanks to the work of J. Arthur \cite{A} and for quasi-split unitary groups thanks to the work of C.-P. Mok \cite{Mk}.)

\null {\bf A.3} Suppose now that $\tau $ is a discrete series representation of a standard Levi subgroup $M$ of $G$. Denote by $\sigma $ a unitary supercuspidal representation of a standard Levi subgroup $M_1$ of $M$ such that, for some $\lambda _1\in a_{M_1}^{M*}$, $\sigma _{\lambda _1}$ is contained in the supercuspidal support of $\tau $ (or, in other words, $\tau $ is a subquotient of $i_{P_1\cap M}^M\sigma_{\lambda _1}$, where $P_1$ is the standard parabolic subgroup with Levi factor $M_1$).

The following lemma follows immediately from the fact that Harish-Chandra's $\mu $-function is defined, up to a constant, as the reciprocal of the scalar operator associated to the composition of the standard intertwining operators with respect to opposite parabolic subgroups \cite{W, V.2} and the product formula \cite{W,V.2.1}:

\null {\bf Lemma:} \it As functions in $\lambda\in a_{M,\Bbb C}^*$, $\mu (\tau _{\lambda })$ and $(\mu /\mu ^M)(\sigma _{\lambda _1+\lambda })$ are directly proportional, i.e. $\mu (\tau _{\cdot })\sim (\mu /\mu ^M)(\sigma _{\lambda _1+\cdot })$. \rm

\null {\bf A.4} Suppose first that there is a supercuspidal representation $\sigma '$ of the corresponding Levi subgroup $M_1'$ of $G'$ in the Vogan $L$-packet of $\sigma $. Then, as the $\mu $-function of a supercuspidal representation is determined, up to a constant, by the reducibility points of the representation, and, by the results of Moeglin \cite{M1, M2} (and of Bernstein-Zelevinsky for $GL_n$) these reducibility points can be read off from the Langlands parameter, one sees that the $\mu $-functions of $\sigma $ and of $\sigma '$ are directly proportional. This implies that, with $P_1'$ equal to the standard parabolic subgroup of $G'$ with Levi factor $M_1'$, $i_{M'\cap P_1'}^{M'}\sigma '_{\varphi(\lambda _1) }$ has a discrete series subquotient $\tau '$, as $\sigma_{\lambda _1}$ is a residue point for $\mu ^G$ (cf. {\bf 2.3}). By the Moeglin-Tadic classification for classical groups and its links to the Langlands correspondence \cite{M1, M2}, one sees that $\tau $ and $\tau '$ have same Langlands parameters, i.e. that they lie in the same Vogan $L$-packet. One concludes by the equality $$\mu(\tau _{\lambda })\sim(\mu/\mu ^M)(\sigma _{\lambda +\lambda _1})\sim(\mu/\mu ^{M'})(\sigma '_{\varphi(\lambda )+\varphi (\lambda _1)})\sim\mu(\tau '_{\varphi (\lambda )}),$$ the first and the last equality following from lemma {\bf A.3}.

\null{\bf A.5} In general, however there may not exist a supercuspidal representation of $M_1'$ in the Vogan $L$-packet of $\sigma $. This happens precisely (see the end of this sub-section), if the generic member in the Vogan $L$-packet is induced from a parabolic subgroup included in the Siegel parabolic subgroup (the one whose standard Levi-factor is isomorphic to a general linear group). To deal with this exceptional case, we will slightly modify the Langlands parameter (and the group $G'$) to have the existence of a supercuspidal representation with good properties and relate this to $M'$.

Before doing that, let us go a little bit in the details of the Langlands correspondence for quasi-split classical groups, which is known thanks to the work of C. Moeglin \cite{M1, M2} on the stablized trace formula for these groups (summarized in \cite{H4, 1.1 and C.3} with the restriction on the literature mentioned in {\bf 3.1}). Denote by $\Cal M_1$ the common $L$-group of $M_1$ and $M'_1$ and by $\eta _1:W_F\times SL_2(\Bbb C)\rightarrow\hbox{\main M}_1$ the Langlands parameter corresponding to the Vogan $L$-packet defined by $\sigma $. Remark that $\Cal M_1$ is a complex reductive group that is a product of general linear groups with a group $\Cal{H}$ of the same type as the $L$-group of $G$. This means that $\eta _1$ is a direct sum of irreducible representations of $W_F$ and a Langlands parameter $\eta _{1,0}:W_F\times SL_2(\Bbb C)\rightarrow\Cal{H}$. There is a canonical inclusion $\iota $ of $\Cal{H}$ into a general linear group. The representation $\iota\circ\eta _{1,0}$ decomposes into irreducible representations of $W_F\times SL_2(\Bbb C)$.

By the restriction of a Langlands parameter to $\Delta W_F$, in notations ${\eta _1}_{\vert \Delta W_F}$ or $(\iota\circ\eta _1)_{\vert \Delta W_F}$, we will mean in the sequel the composition with the map $\Delta:W_F\rightarrow W_F\times SL_2(\Bbb C)$, $\gamma\mapsto (\gamma ,diag(\vert\gamma\vert_F^{1/2}, \vert\gamma \vert_F^{-1/2}))$, while ${\eta _1}_{\vert W_F}$ or $(\iota\circ\eta _1)_{\vert W_F}$ will denote the restriction to the first component.

By the results of Moeglin (see \cite{H4, 1.1, C.3}), it may happen that there is no supercuspidal representation of $M_1'$ in the Vogan $L$-packet of $\sigma $ (by assumption $\sigma $ is a supercuspidal representation of $M_1$, but there may be none on the quasi-split inner form of $M_1$ inside this $L$-packet). In this case, all irreducible components of $(\iota\circ\eta _{1,0})_{\vert W_F}$ have even multiplicities. It can be easily deduced from the above mentioned result of Moeglin that this is equivalent to the property that the supercuspidal support of the $\psi '$-generic member of the Vogan $L$-packet comes from a standard Levi subgroup included in the Siegel parabolic subgroup of $G'$. In fact, the Langlands parameter of the supercuspidal support of the $\psi '$-generic member of the Vogan $L$-packet of $\sigma $ takes the same values as ${\eta_1}_{\vert \Delta W_F}$ and its target space is the smallest Levi subgroup $\Cal{M}_{gen}$ of $\Cal{M}$, containing the image of ${\eta _1}_{\vert \Delta W_F}$.

\null {\bf A.6} To circumvent the problem of not having a supercuspidal representation of $M_1'$ in the Vogan $L$-packet of $\sigma $, we will choose an irreducible representation $\rho $ of $W_F$, which factors through a group of the same type as the $L$-group of $G'$ and which is not a factor of $(\iota\circ\eta _1)_{\vert W_F}$. (As there are infinitely many, this is always possible.) Replacing $\Cal{H}$ by a suitable group $\Cal{H}''$ of bigger rank, but of the same type as $\Cal{H}$, one sees that there is a Langlands parameter $\eta _{1,0}'':W_F\times SL_2(\Bbb C)\rightarrow\Cal{H}''$ such that $\iota\circ\eta _{1,0}''=(\iota\circ\eta_{1,0})\oplus\rho $. Denote by $\Cal{M}_1''$ the complex group obtained from $\Cal{M}_1$, by replacing the factor $\Cal{H}$ by $\Cal{H}''$, and by $G''$ a quasi-split classical group of the same type as $G'$ such that $\Cal{M}_1''$ is the $L$-group of a Levi subgroup of $G''$. Denote by $\eta _1''$ the Langlands parameter $W_F\times SL_2(\Bbb C)\rightarrow\Cal{M}_1''$ obtained from $\eta _1$ by replacing $\eta _{1,0}$ by $\eta _{1,0}''$.

Then $\eta _1''$ defines an $L$-packet of discrete series representations of $M''$ that contains a supercuspidal representation $\sigma ''$. As the non semi-simple factors of $M_1$ and $M_1''$ are the same, there is a canonical isomorphism $\varphi '':a_{M_1,\Bbb C}^*\rightarrow a_{M''_1,\Bbb C}^*$, as for corresponding Levi subgroups of pure inner forms. Denote by $M''$ the standard Levi subgroup of $G''$ that corresponds to $M$ by this identification. As $\rho $ is not a direct summand of $(\iota\circ\eta _1)_{\vert W_F}$, the reducibility points of $\sigma $ and $\sigma ''$ via the isomorphism $\varphi ''$ are the same by the results of Moeglin \cite{M1, M2}. One concludes as in the above case $M''=M'$ (cf. {\bf A.4}) that, with $P_1''$ equal to the standard parabolic subgroup with Levi factor $M_1''$, the representation $i_{P_1''\cap M''}^{M''}(\sigma ''_{\varphi''(\lambda _1)})$ has a discrete series subquotient $\tau ''$ which has, via $\varphi ''$, a $\mu $-function directly proportional to the one of $\tau $. As the $\mu$-function is invariant on (Langlands) $L$-packets of quasi-split classical groups (cf. {\bf A.2}), it follows that the $\mu $-function of $\tau ''$ equals the one of the $\psi '$-generic member $\tau _{gen}''$ in its $L$-packet.

\null {\bf A.7} Denote by $\eta $ the Langlands parameter of $\tau $. Using the relation between the Moeglin-Tadic classification of discrete series representations and the local Langlands correspondence \cite{M1, M2}, one sees that the common Langlands parameter $\eta ''$ of $\tau ''$ and $\tau _{gen}''$ satisfies $\iota\circ\eta ''=(\iota\circ\eta)\oplus\rho $ and that $\rho $ is not a summand of $(\iota\circ\eta)_{\vert W_F}$. Denote by $\tau _{gen}$ the generic member of the Vogan $L$-packet of $\tau $, which is necessarily a representation of $M'$. The values of $\eta _{\vert \Delta W_F}$ (resp. $\eta''_{\vert \Delta W_F}$) are equal to the values of the Langlands parameter of a representation in the supercuspidal support of $\tau _{gen}$ (resp. $\tau ''_{gen}$), and both differ by the summand $\rho $.

Again (as in {\bf A.6, A.4}): as the $\mu $-functions for supercuspidal representations are determined, up to a constant, by the reducibility points of the representations and these points can be read off by the results of Moeglin and Bernstein-Zelevinsky from the Langlands parameter, while $\rho $ is not a direct summand of $(\iota\circ\eta )_{\vert W_F}$, the two generic supercuspidal representations in the supercuspidal support of respectively $\tau_{gen}$ and $\tau ''_{gen}$ have directly proportional $\mu $-functions via the identification $\varphi''\circ\varphi^{-1}$ of $a_{M'}^*$ and $a_{M''}^*$. One concludes from this as in {\bf A.4} by the lemma {\bf A.3} that the $\mu $-functions of $\tau_{gen}$ and $\tau ''_{gen}$ are directly proportional after identification of the underlying spaces.

\null{\bf A.8} We have shown, via the identifications, in {\bf A.6} that the $\mu $-function of $\tau $ is directly proportional to the one of $\tau _{gen} ''$ and in {\bf A.7} that the $\mu $-function of $\tau _{gen}''$ is directly proportional to the one of $\tau _{gen}$. This implies that the $\mu $-functions of $\tau $ and $\tau _{gen}$ are directly proportional.

In conclusion, one sees that, for every discrete series representation $\tau $ of $M$ in the Vogan $L$-packet, there exists a discrete series representation $\tau '$ of $M'$ in the same Vogan $L$-packet with directly proportional $\mu $-function. As, as we recalled above in {\bf A.2}, it is known that the members of a discrete series $L$-packet of a Levi subgroup of a quasi-split classical group have the same $\mu $-function, this finishes the proof of the theorem {\bf A.1}. \hfill{\fin 2}

\enddocument